\documentclass[reqno,11pt]{amsart}
\usepackage{citesort}

\newenvironment{proofl}{\textit{Proof :}}{$\Box$\newline}
\newenvironment{proofp}{\textit{Proof :}}{$\Box$\newline}

\newtheorem{theorem}{Theorem}[section]
\newtheorem{proposition}{Proposition}[section]
\newtheorem{lemma}[theorem]{Lemma}

\theoremstyle{definition}

\theoremstyle{remark}
\newtheorem{remark}[theorem]{Remark}

\numberwithin{equation}{section}

\def\R{{\mathbb{R}}} 
\def\N{{\mathbb{N}}}

\def\epsilon{{\varepsilon}}
\def\phi{{\varphi}}
\def\theta{{\vartheta}}
\def\scal{{H_2}}
\def\scalk{{H_k}}
\def\ov{{\overline{v}}}
\def\ou{{\overline{u}}}
\def\oW{{\overline{W}}}
\def\Omb{{\overline{\Omega}}}
\def\H{\stackrel{\sim}{H}}
\def\F{{F}}
\def\n{\stackrel{\rightarrow}{n}}

\begin{document}
\title{Entire spacelike hypersurfaces of prescribed scalar curvature in Minkowski space}

\author{Pierre Bayard}
\address{Pierre Bayard: Instituto de F\'{\i}sica y Matem\'aticas. U.M.S.N.H. Ciudad Universitaria. CP. 58040 Morelia, Michoacan, Mexico.}
\email{bayard@ifm.umich.mx}


\date{\today.}


\begin{abstract}
We prove existence and uniqueness of entire spacelike hypersurfaces in the Minkowski space with prescribed negative scalar curvature, and with given values at infinity which stay at a bounded distance of a lightcone.
\end{abstract}
\maketitle
\markboth{PIERRE BAYARD}{ENTIRE HYPERSURFACES OF PRESCRIBED SCALAR CURVATURE}
\textit{Mathematical Subject Classification (2000):} 35J60, 53C50.
\section{Introduction}
In this article, we study existence and uniqueness in the Minkowski space of entire hypersurfaces of prescribed scalar curvature and with given values at infinity. The Minkowski space $\R^{n,1}$ is the space $\R^n\times\R$ with the metric$$ds^2=dx_1^2+\cdots+dx_n^2-dx_{n+1}^2.$$
An entire spacelike hypersurface $M$ is a graph of a smooth function $u$ defined over $\R^n$ such that $|Du|<1$ on $\R^n.$ The latter means equivalently that the induced metric on $M$ is riemannian. In the coordinates $x_1,\ldots,x_n,$ the metric on $M$ is given by the matrix
$$\left(\delta_{ij}-u_iu_j\right)_{ij},$$
and the curvature endomorphism by
\begin{equation}\label{endcurv}
\frac{1}{\sqrt{1-|Du|^2}}\left(\delta_{ij}+\frac{u_iu_j}{1-|Du|^2}\right)_{ij}D^2u.
\end{equation}
By definition, the principal curvatures $\lambda_1,\ldots,\lambda_n$ are the eigenvalues of the curvature endomorphism (\ref{endcurv}), and the $k^{th}$ mean curvature is the $k^{th}$ elementary symmetric function of the principal curvatures, and is denoted by:
$$\scalk[u]=\sigma_k(\lambda_1,\ldots,\lambda_n)=\sigma_k.$$
Let us first mention articles concerning entire solutions of the prescribed mean curvature equation $(k=1)$ in lorentzian manifolds : in the Minkowski space, Cheng and Yau \cite{CY} proved the Bernstein property for the entire maximal hypersurfaces and Treibergs \cite{T} obtained the classification of the entire constant mean curvature hypersurfaces; Bartnik \cite{B} proved the existence of maximal hypersurfaces in asymptotically flat spacetimes and also obtained estimates for the Dirichlet problem for the prescribed mean curvature equation.

In \cite{De}, Ph. Delano\"e solved the Dirichlet problem for the prescribed Gauss curvature equation $(k=n)$ in Minkowski space when the datas are strictly convex. In \cite{G}, Bo Guan solved this problem under the weaker assumption that there is a strictly convex subsolution. He also solved the Plateau problem. O. Schn\"urer \cite{Sch} solved the Dirichlet problem in lorentzian manifolds. 

Ph. Delano\"e mentioned to us the problem of finding entire hypersurfaces of prescribed Gauss curvature in Minkowski space asymptotic to a lightcone, and O. Schn\"urer explained to us how to obtain the uniform interior gradient estimate required to solve this problem (unpublished; see remark \ref{rmk2} below). We are grateful to the Referee who mentioned to us An-Min Li article \cite{L}, where the author proved the existence of entire convex spacelike hypersurfaces of prescribed positive Gauss curvature and with bounded principal curvatures. He also proved the completeness of these hypersurfaces equipped with the induced metric. 

We are interested here in the scalar curvature $S$ of the spacelike hypersurface $M,$ which is linked to $\sigma_2$ by : 
$$S=-2\sigma_2.$$
Our approach of the existence of entire hypersurfaces of prescribed scalar curvature relies in a crucial way on previous works concerning the Dirichlet problem: we solved in \cite{Bay} the Dirichlet problem in the ambient Minkowski space of dimension 4; J.Urbas \cite{U} proved a maximum principle for the curvature of solutions, which permits him to solve the Dirichlet problem in all dimensions. He also obtained an interior curvature bound which implies the existence of locally smooth solutions, in the case of spacelike affine boundary data. The articles quoted in this introduction, especially \cite{B},\cite{T},\cite{U}, are at the source of the present work. For further references concerning the large literature of problems of prescribed curvature in lorentzian manifolds, please see references in \cite{Bay} or in Schn\"urer's article \cite{Sch}.

We will say that the function $u$ is admissible if $u$ is spacelike and the principal curvatures of its graph belong to
$$\Gamma_2=\{\lambda\in\R^n\ s.t.\ \sigma_1(\lambda)>0\mbox{ and } \sigma_2(\lambda)>0\}.$$  
It is well known that $\scal$ defines a fully non linear elliptic operator such that $\scal^{\frac{1}{2}}$ is a concave function of $D^2u,$ if $u$ is admissible.
\\
\\Our first result concerns existence and uniqueness of entire hypersurfaces of prescribed scalar curvature asymptotic to a lightcone. We assume the existence of lower and upper barriers which are asymptotic to the lightcone. 
\begin{theorem}\label{the0}
Let $H$ be a positive function $C^{2,\alpha}$ for some $\alpha\in\ ]0,1[$ on the upper cone $C^+=\{(x,x_{n+1})|\ x_{n+1}>|x|\}.$ Let us assume that there exist $\varphi_1\in\ C^{4,\alpha}(\R^n),$ strictly convex and spacelike, and $\varphi_2\in\ C^2(\R^n),$ spacelike, such that
$$\scal[\varphi_1]\geq H(.,\varphi_1),\ \scal[\varphi_2]\leq H(.,\varphi_2)\mbox{ in }\R^n,$$
$$\lim_{|x|\rightarrow +\infty}\varphi_i(x)-|x|=0 \mbox{ for }i=1,2,$$
and $\varphi_1\leq\varphi_2.$ Then there exists a function $u:\R^n\rightarrow\R,$ belonging to $C^{4,\alpha},$ spacelike, such that 
\begin{equation}\label{eqcurv}
\scal[u]=H(.,u)\mbox{ in }\R^n,
\end{equation}
\begin{equation}\label{condinfty}
\lim_{|x|\rightarrow +\infty}u(x)-|x|=0,
\end{equation}
and
\begin{equation}\label{encadr}
\varphi_1\leq u\leq \varphi_2\mbox{ in }\R^n.
\end{equation}
Moreover, if $\frac{\partial H}{\partial x_{n+1}}\geq 0,$ the solution of (\ref{eqcurv})-(\ref{condinfty}) is unique and satisfy (\ref{encadr}).
\end{theorem}
\begin{remark}
Since the barriers $\varphi_1$ and $\varphi_2$ are supposed to be spacelike and asymptotic to the lightcone, their graphs lie in the upper cone $C^+.$ 
\end{remark}
\begin{remark}\label{rmk1}
If we have a pinching $h_1\geq H\geq h_2$ where $h_1,h_2$ are positive constants, setting $ \alpha_1=\sqrt{\frac{n(n-1)}{2h_1}}$ and  $\alpha_2=\sqrt{\frac{n(n-1)}{2h_2}},$ the functions $\varphi_1,\varphi_2$ defined by $\varphi_1(x)=\sqrt{\alpha_1^2+|x|^2},$ and $\varphi_2(x)=\sqrt{\alpha_2^2+|x|^2}$ are natural barriers: the graph of $\varphi_1$ (resp. $\varphi_2$) is the upward hyperbolo\"{\i}d asymptotic to the lightcone $x_{n+1}=|x|$ whose $\sigma_2$ is $h_1$ (resp. $h_2$).
\end{remark}
Our second result concerns the existence and uniqueness of entire hypesurfaces of prescribed scalar curvature asymptotic to a small perturbation of a lightcone. We assume here a pinching condition on the prescribed curvature.
\begin{theorem}\label{the1}
Let $H\in C^{2,\alpha}(\R^n\times\R),$ $\alpha\in\ ]0,1[,$ such that $h_1\geq H\geq h_2$ for some positive constants $h_1,h_2.$ Let $f$ be a function of class $C^2$ on the unit sphere of $\R^n.$ Then  there exists a function $u:\R^n\rightarrow\R,$ belonging to $C^{4,\alpha},$ spacelike, such that 
\begin{equation}\label{eqn1}
\scal[u]=H(.,u)\mbox{ in }\R^n,
\end{equation}
and
\begin{equation}\label{eqn2}
\lim_{|x|\rightarrow +\infty}u(x)-|x|-f\left(\frac{x}{|x|}\right)=0.
\end{equation}
Moreover, if $\frac{\partial H}{\partial x_{n+1}}\geq 0,$ the solution of (\ref{eqn1})-(\ref{eqn2}) is unique.
\end{theorem}
\begin{remark}
Theorem \ref{the1} implies in particular that there is a lot of entire spacelike hypersurfaces of constant negative scalar curvature in Minkowski space, with projective boundary values at infinity given by a lightcone. For the notion of projective boundary values at infinity, and the statement in the mean curvature case, see Treibergs \cite{T} paragraph 6. The classification of entire hypersurfaces with constant scalar curvature by their values at infinity, as in \cite{T} for the mean curvature case, seems an interesting open question.  
\end{remark}
\begin{remark}
We thank the Referee for the following remark : it will be clear from the proof of theorem \ref{the0} that if there exist barriers $\varphi_1$ and $\varphi_2$ as in the statement of theorem \ref{the0}, but with
\begin{equation}\label{asymptalpha}
\lim_{|x|\rightarrow +\infty}\left(\varphi_i(x)-\alpha|x|\right)=0 \mbox{ for }i=1,2,
\end{equation}
for some $\alpha\in\ ]0,1[,$ then there is an admissible solution of (\ref{eqcurv})-(\ref{condinfty}) satisfying (\ref{encadr}). Such a generalization of theorem \ref{the1} is not possible because the asymptotic behaviour (\ref{asymptalpha}) is not compatible with the pinching assumption on $H$ unless $\alpha=1$ (for this last claim, see for example Treibergs \cite{T} section 6, step 1 of the proof of theorem 1).
\end{remark}

\begin{remark}
Let us conclude with some easy observations concerning the solutions of (\ref{eqcurv})-(\ref{condinfty}) or of (\ref{eqn1})-(\ref{eqn2}).
\\{\bf 1.} A solution $u$ is necessarily admissible : let us first recall the useful Mac-Laurin's inequality : if $\sigma_1,\sigma_2\geq 0,$ we have :
\begin{equation}\label{inegml}
\sigma_2\leq \frac{n-1}{2n}\sigma_1^2.
\end{equation}
Since $u\rightarrow_{|x|\rightarrow+\infty} +\infty,$ $u$ attains its minimum at some point. At this point, since $\scal[u]>0$ and in view of the Mac-Laurin inequality (\ref{inegml}), $u$ is admissible. The set of points where $u$ is admissible is an open set, and is closed (using $\scal[u]>0$ together with the Mac-Laurin inequality). It is thus $\R^n.$ 
\\{\bf 2.} Let us  prove uniqueness under the assumption $\frac{\partial H}{\partial x_{n+1}}\geq 0,$ by contradiction : let $u$ and $v$ be two (admissible) solutions. If there exists $x_0$ such that $u(x_0)<v(x_0),$ let $\delta>0$ be such that $\delta<v(x_0)-u(x_0),$ and set $u_{\delta}=u+\delta.$ Since $\lim_{|x|\rightarrow +\infty} u_{\delta}(x)-v(x)=\delta,$ we see that $\Omega_{\delta}=\{u_{\delta}<v\}$ is a bounded open set. We have $v=u_{\delta}$ on $\partial\Omega_{\delta},$ $\scal[v]=H(.,v)$ and $\scal[u_{\delta}]=\scal[u]=H(.,u)\leq H(.,u_{\delta})$ on $\Omega_{\delta}$ since $\frac{\partial H}{\partial x_{n+1}}\geq 0.$ By a standard comparison principle ($v$ admissible and $\frac{\partial H}{\partial x_{n+1}}\geq 0$), we thus obtain $v\leq u_{\delta}$ on $\Omega_{\delta},$ which is impossible.
\\{\bf 3.} From elliptic regularity theory, if $H\in\ C^{\infty}(\R^n\times\R),$ a solution $u$ belongs to $C^{\infty}(\R^n).$
\end{remark}
This article is organized as follows : we first focus on the proof of theorem \ref{the0}. In section \ref{sec00}, we solve the Dirichlet problem for the prescribed scalar curvature equation between two barriers. In section \ref{sec0}, we present the method of construction of an entire solution, putting forward the required a priori estimates. We obtain the gradient estimate in section \ref{gradient} (our main contribution), and the $C^2$ estimate in section \ref{c2}. The theorem \ref{the1} will easily follow from these proofs once barriers and auxiliary functions are constructed. The construction of such functions is the aim of section \ref{lastsec}.

\section{The Dirichlet problem with barriers}\label{sec00}

We know from \cite{Bay},\cite{U} that we can solve the Dirichlet problem for the prescribed scalar curvature equation in the Minkowski space if the open set $\Omega$ and the boundary data $\varphi$ are uniformly convex, and if the additional condition $\frac{\partial H}{\partial x_{n+1}}\geq 0$ holds . Our purpose here is to remove the latter condition, and to prove that if there exist lower and upper barriers we can find a solution of the Dirichlet problem which lies between these barriers. Using theorem 1.1 of Urbas \cite{U} and a new fixed point argument indicated to us by Ph. Delano\"e \cite{De2}, we prove the following.
\begin{theorem}\label{theo3}
Let $\Omega$ be a uniformly convex domain in $\R^n$ with $\partial\Omega$ $C^{4,\alpha}$ for some $\alpha\in\ ]0,1[,$ and $H\in\ C^{2,\alpha}(\overline{\Omega}\times\R)$ be a positive function. Let $\varphi_1\in\ C^{4,\alpha}(\Omb),$ strictly convex and spacelike, and $\varphi_2\in\ C^2(\Omb),$ spacelike, such that 
$$\scal[\varphi_1]\geq H(.,\varphi_1),\ \scal[\varphi_2]\leq H(.,\varphi_2)\mbox{ in }\Omega,$$
and $\varphi_1<\varphi_2$ in $\Omb.$ Then there exists a spacelike function $u$ belonging to $C^{4,\alpha}({\overline{\Omega}})$ such that  
\begin{equation}\label{pd}
\left\{\begin{array}{l}
\scal[u]=H(.,u)\mbox{ in }\Omega,\\
u=\varphi_1\mbox{ on }\partial \Omega,
\end{array}\right.
\end{equation}
and $\varphi_1\leq u\leq\varphi_2.$
\end{theorem}
Before the proof of the theorem, let us recall the following comparison principles for the operator $\scal.$ We omit the (classical) proofs. 
\begin{proposition}\label{comppple}
Let $\Omega$ be a bounded domain with $\partial\Omega\in\ C^2,$ $H$ be a function belonging to $C^1(\Omb\times\R)$ such that $\frac{\partial H}{\partial x_{n+1}}\geq 0,$ and let $u,v\in\ C^2(\Omb)$ be two spacelike functions. Let us assume that $v$ is admissible, $\scal[v]\geq H(.,v),$ $\scal[u]\leq H(.,u)$ in $\Omega,$ and $v\leq u$ on $\partial\Omega.$ Then 
$$v\leq u \mbox{ in }\Omb.$$
Moreover, if $v(x_0)=u(x_0)$ at an interior point $x_0\in\Omega,$ then $v\equiv u$ in $\Omb$ (the strong comparison principle), and if $v(x_0)=u(x_0)$ at $x_0\in\ \partial\Omega$ and $v<u$ in $\Omega$ then $(u-v)_{\n}(x_0)>0$ (the Hopf lemma).
\end{proposition} 
In the proposition and below, we write the index $\n$ to denote the interior normal derivative at a boundary point.
\\
\\\textit{Proof of theorem \ref{theo3} :} We suppose that $\varphi_1$ is not a solution of (\ref{pd}). We denote by $K$ the compact set
$$K=\{(x,z)|\ x\in\ {\overline{\Omega}},\ \varphi_1(x)\leq z \leq\varphi_2(x)\}.$$ 
Setting 
$$k=\max\left(\sup_K\frac{1}{H}\frac{\partial H}{\partial x_{n+1}},0\right),$$
the function $z\mapsto H(x,z)e^{-kz}$ is decreasing on $[\varphi_1(x),\varphi_2(x)] $ for all $x\in\ {\overline{\Omega}}.$
\\The usual iterative method (see \cite{CH} p369-371) fails here, and a new method is necessary (aim of \cite{De2}).
\\Let us consider the Banach space
$$E=\{\ov\in\ C^{2,\alpha}(\Omb)|\ \ov=0\mbox{ on }\partial\Omega\}$$
and the convex open set of $E$
$$W=\{\ov\in\ E|\ \ov>0\mbox{ in } \Omega,\ \ov_{\n}>0\mbox{ on }\partial\Omega,\ \mbox{and }\ov<\varphi_2-\varphi_1\mbox{ on }\Omb\}.$$
We then define the operator $\begin{array}[t]{rcl}
T:[0,1]\times\oW & \rightarrow & E\\
(t,\ov) & \mapsto & \ou,
\end{array}$
where $\ou\in\ E$ is such that $u=\ou+\varphi_1$ is the admissible solution (belonging to $C^{4,\alpha}$) of the Dirichlet problem
\begin{eqnarray}
\scal[u]e^{-ku}&=&tH(.,v)e^{-kv}+(1-t)H(.,\varphi_1)e^{-k\varphi_1}\mbox{ in }\Omega\label{e1}\\
u&=&\varphi_1\mbox{ on }\partial\Omega.\label{e2}
\end{eqnarray}
Here $v=\ov+\varphi_1.$ This operator $T$ is well defined since the positive sign of the right hand side term of (\ref{e1}) and the non-negativity of $k$ allow to apply theorem 1.1 of Urbas \cite{U}. The aim is now to prove that $T(1,.)$ has a fixed point. We first prove that $T$ takes its values in the set $W,$ we then estimate the fixed points of $T,$ and we conclude thanks to the fixed point theorem of Browder-Potter \cite{BP}. 
\\1- \textit{$T$ takes its values in the set $W :$} we note that $\varphi_1$ and $\varphi_2$ are respectively sub and supersolution of (\ref{e1})-(\ref{e2}) since
$$\scal[\varphi_1]e^{-k\varphi_1}\geq H(.,\varphi_1)e^{-k\varphi_1}\geq H(.,v)e^{-kv}\geq H(.,\varphi_2)e^{-k\varphi_2}\geq \scal[\varphi_2]e^{-k\varphi_2}$$
$(\varphi_2\geq v\geq\varphi_1$ and $z\mapsto H(.,z)e^{-kz}$ is decreasing). The comparison principle (proposition \ref{comppple}) thus implies that $\varphi_1\leq u\leq\varphi_2$ in $\Omega.$ We remark that these inequalities are strict : if $u=\varphi_1$ at an interior point then, from the strong comparison principle, $u \equiv\varphi_1$ and thus :
$$\scal[\varphi_1]e^{-k\varphi_1}=tH(x,v)e^{-kv}+(1-t)H(x,\varphi_1)e^{-k\varphi_1}\leq H(.,\varphi_1)e^{-k\varphi_1}$$
since $\varphi_1\leq v\leq\varphi_2$ and $z\mapsto H(.,z)e^{-kz}$ decreases. Thus $\scal[\varphi_1]\leq H(.,\varphi_1)$ and $\varphi_1$ is a solution, which is impossible. Thus $\varphi_1<u$ in $\Omega.$ If $u=\varphi_2$ at an interior point of $\Omega,$ then $u\equiv \varphi_2$ in $\Omega,$ which is impossible since $\varphi_1<\varphi_2$ on $\partial\Omega.$ We finally need to prove that $\ou_{\n}>0$ on $\partial\Omega :$ we know that $\scal[\varphi_1]e^{-k\varphi_1}\geq\scal[u]e^{-ku}$ in $\Omega,$ and $u=\varphi_1$ on $\partial\Omega.$ Since $u>\varphi_1$ in $\Omega,$ the Hopf lemma implies that $(u-\varphi_1)_{\n}>0$ on $\partial\Omega,$ which achieves the proof.
\\2- \textit{The fixed points of $T$ are under control :} the fixed points of $T$ are solutions of the Dirichlet problem :
\begin{eqnarray}
\scal[u]&=&\H(.,u)\mbox{ in }\Omega\label{fix1}\\
u&=&\varphi_1\mbox{ on }\partial\Omega,\label{fix2}
\end{eqnarray}
where
\begin{equation}\label{Htilde}
\H(.,u)=tH(.,u)+(1-t)H(.,\varphi_1)e^{k(u-\varphi_1)}.
\end{equation}
Let us examine the $C^1$ and the $C^2$ estimates obtained in \cite{Bay} and in \cite{U}, specifying their dependence on $\H$ and its derivatives.
\\\textit{The $C^1$ estimate :} the maximum principle and the barriers construction at the boundary obtained in \cite{Bay} section 3 readily extend to the case where the function of prescribed curvature $\H$ also depends on $u.$ We easily obtain :
\begin{equation}\label{C1pd}
\sup_{\Omb}|Du|\leq 1-\theta,
\end{equation} 
where $\theta=\theta(\inf_K\H,\|\H\|_{1,K},\sup_{\Omb}(\varphi_2-\varphi_1)).$
\\
\\\textit{The $C^2$ estimate :} the maximum principle of Urbas \cite{U} theorem 1.2 embraces the case of equations (\ref{fix1})-(\ref{fix2}) and gives the estimate :
\begin{equation}\label{C2-1pd}
\sup_{\Omb}|D^2u|\leq C_1,
\end{equation}
where $C_1$ depends on $\sup_{\Omb}|Du|,$ $\sup_{\partial\Omega}|D^2u|,$ $\|\H\|_{2,K},$ $\|\varphi_1\|_{1,\Omb},$ and on a positive lower bound on the minimum eigenvalue of $D^2\varphi_1$ on $\Omb.$
\\We now study the barriers constructions at the boundary given in \cite{Bay} section 4-2 : with the key lemma 4.2 at hand (cf. \cite{Bay} page 17), the barriers constructions proposed in \cite{Bay} to estimate the mixed second derivatives and the normal second derivative at the boundary remain unchanged. Following the notations adopted there, we write the equation of prescribed scalar curvature in the form ${\F}_2[u]=f(.,u,Du),$ where $f(.,u,Du)=(1-|Du|^2)H(.,u).$ The only new fact here is that the function $f$ also depends on $u$ (and not only on $x$ and $Du$). We used in the proof of the key inequality (26) of \cite{Bay} the equation of prescribed curvature differentiated once (formula (25)) to substitute third derivatives of $u.$ Using the same device here, the new terms that we get only depend on the first derivatives of $u,$ and are thus obviously under control.This is the only required modification in the proof, and we thus get from \cite{Bay} the estimate:
\begin{equation}\label{C2-2pd}
\sup_{\partial\Omega}|D^2u|\leq C_2,
\end{equation}
where $C_2$ depends on $\sup_{\Omb}|Du|,$ $\|\H\|_{1,K},$ $\|\varphi_1\|_{4,\Omb}.$
\\The estimates (\ref{C1pd}),(\ref{C2-1pd}),(\ref{C2-2pd}) and the expression (\ref{Htilde}) of $\H$ imply a $C^2$ estimate for the fixed points of $T$ which is independent of the parameter $t\in\ [0,1].$ The Evans-Krylov and the Schauder theories yield a $C^{4,\beta}$ estimate; in particular, we get the existence of a constant $C$ such that : for all $\ou\in\ W$ fixed point of $T,$
$$\|\ou\|_{2,\alpha}<C.$$
The constant $C$ depends on $\sup_{\Omb}(\varphi_2-\varphi_1),\inf_KH,\|H\|_{2,K},\|\varphi_1\|_{4,\Omb},$ and on a positive lower bound on the minimum eigenvalue of $D^2\varphi_1$ on $\Omb.$
\\3- \textit{$T(1,.)$ has a fixed point :} let us set $W_c=\{\ov\in\ W|\ \|\ov\|_{2,\alpha}<C\},$ and consider $T:[0,1]\times\oW_c\rightarrow E.$ The set $W_c$ is convex.
\\a- $T$ is a continuous and compact map in view of the above estimates on the solutions of the Dirichlet problem (\ref{e1})-(\ref{e2}).
\\b- $T(0,.)$ maps $\partial W_c$ in $\oW_c:$ $T(0,.)\equiv \ou_0$ where $u_0=\ou_0+\varphi_1$ is the unique admissible solution of
\begin{eqnarray}
\scal[u]&=&H(.,\varphi_1)e^{k(u-\varphi_1)}\mbox{ in }\Omega,\\
u&=&\varphi_1\mbox{ on }\partial\Omega.
\end{eqnarray}
The function $\ou_0$ belongs to $W$ since we proved that $ImT\subset W,$ and is thus a fixed point of $T(0,.).$ Thus $\|\ou_0\|_{2,\alpha}<C,$ and $\ou_0\in\ W_c.$ 
\\c- $T(t,.)$ does not have any fixed point on $\partial W_c,$ since, as seen above, any fixed point of $T(t,.)$ belongs to the open set $W_c$ by construction of $C.$
\\The theorem of Browder-Potter \cite{BP} implies that $T(1,.)$ has a fixed point (belonging to $C^{4,\alpha}$), which achieves the proof of the theorem.
$\Box$
\newline

\section{The construction of an entire solution}\label{sec0}

We present here the principle of our construction of an entire solution of (\ref{eqcurv})-(\ref{condinfty}). We suppose that the barriers $\varphi_1$ and $\varphi_2$ are not solutions. 
\\Let us first prove that this implies that $\varphi_1<\varphi_2$ in $\R^n :$ indeed, by contradiction, suppose that $\varphi_1(x_0)=\varphi_2(x_0)$ for some point $x_0\in\R^n,$ and let us fix a large ball $B$ containing $x_0,$ and a constant $k$ sufficiently large such that $z\mapsto H(x,z)e^{-kz}$ is decreasing on $[\varphi_1(x),\varphi_2(x)] $ for all $x\in\ {\overline{B}}$ (as in the proof of theorem \ref{theo3}). We have
$$\scal[\varphi_1]e^{-k\varphi_1}\geq H(.,\varphi_1)e^{-k\varphi_1}\geq H(.,\varphi_2)e^{-k\varphi_2}\geq \scal[\varphi_2]e^{-k\varphi_2}$$
on $\overline{B}.$ Since $\varphi_1\leq\varphi_2$ on $\partial B,$ the strong comparison principle implies that $\varphi_1\equiv\varphi_2$ on $\overline{B}.$ Since the ball $B$ is arbitrarily large, we conclude that $\varphi_1\equiv\varphi_2$ on $\R^n,$ and the barriers would thus be solutions of (\ref{eqcurv})-(\ref{condinfty}), which is excluded.
\\
\\For any positive $R,$ denoting by $B_R$ the ball in $\R^n$ centered at 0 and with radius $R,$ we set $u_R$ for an admissible solution of  
$$\left\{\begin{array}{l}
\scal[u_R]=H(.,u_R)\mbox{ in }B_R,\\
u_R=\varphi_1\mbox{ on }\partial B_R
\end{array}\right.$$
such that $\varphi_1\leq u_R\leq\varphi_2$ (theorem \ref{theo3}). We need to prove the following locally uniform estimates :
for any $R_0\geq 0,$ there exist $R_1=R_1(R_0)$ sufficiently large, $\theta\in\ ]0,1[,$ and $C\geq 0$ such that : for every $R\geq R_1,$
\begin{equation}
\sup_{B_{R_0}}|Du_R|\leq 1-\theta,\mbox{ and }\|u_R\|_{2,B_{R_0}}\leq C.
\end{equation}
With these estimates at hand, Evans-Krylov interior second derivative H\"older estimate, and Schauder interior regularity theory imply locally uniform $C^{4,\alpha}$ estimates. A diagonal process then yields a subsequence $u_{R_k},$ $R_k\rightarrow+\infty,$ that locally converges to a solution $C^{4,\alpha}$ of (\ref{eqcurv})-(\ref{condinfty}).
\section{The uniform gradient estimate}\label{gradient}
We begin by defining some notations. The standard basis of $\R^{n,1}$ will be denoted by $\{e_1,\ldots,e_n,e_{n+1}\}.$ We denote the metric paring by $<.,.>,$ and the associated norm of spacelike vectors by $|.|.$ Let $M$ be a spacelike hypersurface and $N$ be its future-directed unit normal. The second fundamental form $II$ of $M$ is thus given by : $\forall p\ \in\ M,$ $\forall X,Y\ \in\ T_pM,$
\begin{equation}\label{ffal}
II_p(X,Y)=<X,D_YN>,
\end{equation}
where $D$ is the connection in $\R^{n,1}.$ If $(\hat{e_1},\ldots,\hat{e_n}),$ is a frame on $M,$ we will denote the components of $II$ in this frame by $h_{ij}.$ We will raise or lower indices with respect to the induced metric on $M.$ The components of the endomorphism of curvature are thus denoted by $h^i_j,$ and we will write the equation of prescribed scalar curvature in the form :
\begin{equation}\label{eqintr}
F\left((h^i_j)_{i,j}\right)=H,
\end{equation}
where $F(A)$ is the sum of the principal minors of order 2 of the matrix $A.$ Let
$$F_i^j=\frac{\partial F}{\partial h^i_j}\left((h^i_j)_{i,j}\right).$$
If $(h^i_j)_{i,j}$ is diagonal, so is $(F^j_i)_{i,j},$ and
$$ (F^j_i)_{i,j}=Diag(\sigma_{1,1},\sigma_{1,2},\ldots,\sigma_{1,n}),$$
where $\forall i,$ $\displaystyle{\sigma_{1,i}=\sum_{k,k\neq i}\lambda_k.}$ Moreover $(F^j_i)_{i,j}$ defines a (1,1) tensor on $M.$ Raising the index $i$ thanks to the induced metric on $M,$ we also will use the symmetric tensor $(F^{ij})_{i,j}.$
\\
\\We denote the induced connection on $M$ by $\nabla.$ If $f$ is a smooth function defined on a neighborhood in $\R^{n,1}$  of $M,$ $grad_Mf=\nabla f$ is given by :
$$\nabla f=Df+<Df,N>N,$$
and, if $M$ is the graph of a smooth function $u,$ the vector $N$ is given by :
$$N=\nu (Du,1),$$
where
\begin{equation}\label{nu}
\nu=-<e_{n+1},N>=\frac{1}{\sqrt{1-|Du|^2}}.
\end{equation}
We thus easily get the estimate :
$$|\nabla f|\leq\nu |Df|_{eucl},$$
where $|.|_{eucl}$ denotes the euclidean norm in $\R^{n+1}.$ For $u$ itself, we get the identities :
$$|\nabla u|=\nu |Du|=\sqrt{\nu^2-1}.$$
Let us recall the useful Gauss formula : for all $\alpha=1,\ldots,n+1,$
\begin{equation}\label{gauss1}
\nabla^2(x_{\alpha})=N_{\alpha}II,
\end{equation}
where $x_{\alpha}$ is the $\alpha^{th}$-coordinate function, and $N_{\alpha}$ the $\alpha^{th}$-component of $N.$
\\For $\alpha=n+1,$ since ${x_{n+1}}_{|M}$ is the function $u,$ we get : 
\begin{equation}\label{gauss2}
\nabla^2u=\nu II.
\end{equation}
We will denote the components of covariant derivatives with a semi-colon, or we will omit it, if there is no risk of confusion. We will for instance denote $\nabla_{\hat{e_k}}II\left(\hat{e_i},\hat{e_j}\right)$ by $h_{ij;k},$ and, if $\eta$ is a function on $M,$ $\nabla^2_{\hat{e_i},\hat{e_j}}\eta$ by $\eta_{ij}.$ 
\\
\\The aim of this section is to prove the following locally uniform gradient estimate : 
\begin{proposition} 
 Given any radius $R_0,$ there exists $\theta_{R_0}\in\ ]0,1]$ and $R_1=R_1(R_0)$ such that : $\forall R\geq R_1,$ 
\begin{equation}\label{estgrad}
\sup_{B_{R_0}}|Du_R|\leq 1-\theta_{R_0}.
\end{equation}
\end{proposition}
Our estimate relies on the existence of a spacelike function $\psi$ satisfying
\begin{equation}\label{ptepsi}
\psi\leq\varphi_1-\delta\mbox{ on }B_{R_0}\mbox{ and }\psi>\varphi_2\mbox{ on }\R^n\backslash B_{R_1},
\end{equation}
where $\delta$ is a controlled positive constant, and $R_1>R_0.$ We may construct $\psi$ as follows : setting $\delta_0=\inf_{B_{R_0}}(\varphi_1(x)-|x|)$ and $\delta=\frac{\delta_0}{2},$ the function $\psi=\varepsilon+\sqrt{\varepsilon^2+|x|^2}$ where $\varepsilon=\frac{\delta}{2}$ is such that
$$\psi\leq\varphi_1-\delta\mbox{ on }B_{R_0}.$$
Since we have 
$$\lim_{|x|\rightarrow+\infty}\psi(x)-|x|=\varepsilon,$$
there exists $R_1=R_1(R_0,\varphi_1,\varphi_2)$ such that 
$$\psi>\varphi_2\mbox{ on }\R^n\backslash B_{R_1}.$$
When taking below extrema of continuous functions of $(x,t)\in\ \R^{n+1}$ on $\{\psi\leq\varphi_2\},$ we will mean extrema taken on the {\textit compact} set 
$$\{(x,t),\ x\in\ B_{R_1},\ \psi(x)\leq t\leq \varphi_2(x)\}.$$
Let $R\geq R_1,$ and let us denote $u_R$ by $u.$ Let $\Gamma$ be the piece of the graph $M$ of $u$ on which $u\geq \psi.$ Following ideas of \cite{B}, we set $\varphi=\eta\nu$ on $\Gamma,$ with 
\begin{equation}\label{defeta}
\eta=(e^u-e^\psi)^K,
\end{equation}
where $K$ is a large constant to be chosen, and $\nu$ is defined by (\ref{nu}).

The function $\varphi$ is non-negative, null on $\partial \Gamma,$ and thus attains its maximum at an interior point $p$ of $\Gamma.$ We prove the following lemma :
\begin{lemma}\label{lem}
Let $\alpha>0$ such that $1-\alpha>\sup_{\psi\leq\varphi_2}|D\psi|.$ For $K$ sufficiently large, we have at $p,$ maximum interior of $\varphi$ :
\begin{equation}\label{estnu}
\nu\leq\frac{1}{\sqrt{1-\left(\frac{\sup_{\psi\leq\varphi_2}|D\psi|}{1-\alpha}\right)^2}}.
\end{equation}
The constant $K$ depends on $\alpha,$ $\displaystyle{\sup_{\psi\leq\varphi_2}(\varphi_2-\psi),}$ $\displaystyle{\inf_{\psi\leq\varphi_2}H}$ and $\displaystyle{\sup_{\psi\leq\varphi_2}|DH|_{eucl}}.$
\end{lemma}
This lemma implies the estimate (\ref{estgrad}), since we have $\eta\leq C_1$ on the set $\{\psi\leq\varphi_2\},$ and, in view of (\ref{ptepsi}), $\eta\geq C_2$ on $B_{R_0},$ where $C_1,C_2$ are positive constants under control.
\\
\\To prove lemma \ref{lem}, we first write the conditions expressing that $\varphi$ is maximum at $p$ (lemma \ref{lem0}, lemma \ref{lem5} inequality (\ref{star4})). We then estimate the terms in (\ref{star4}) : it is the aim of the lemmas \ref{lem1} and \ref{lem4}. These lemmas lead to the key inequality (\ref{star2}) which permits to achieve the proof.
\\
\\We choose at $p$ an orthonormal basis of principal directions $(\hat{e_1},\ldots,\hat{e_n})$ with ordered principal curvatures $\lambda_1\geq\lambda_2\geq\cdots\geq\lambda_n.$
\begin{lemma}\label{lem0}
At $p,$ we have : $\forall i,$
\begin{equation}\label{extr3}
\eta u_i\lambda_i+\eta_i\nu=0.
\end{equation}
\end{lemma}
\begin{proofl}
Since $\varphi$ reaches its maximum at the interior point $p,$ we have : $\forall i,$ $\varphi_i=0,$ and thus,
\begin{equation}\label{extr1}
\nu_i=-\frac{\eta_i}{\eta}\nu.
\end{equation}
 Denoting by $T$ the component tangential to $M$ of $e_{n+1},$ we have, from (\ref{ffal}) and (\ref{nu}), $d\nu=-II(T,.).$ Moreover, since $u={x_{n+1}}_{|M},$ and since the minkowskian gradient of $x_{n+1}$ is $-e_{n+1},$ we have $\nabla u=-T.$
\\Hence 
\begin{equation}\label{nui}
\nu_i=u_i\lambda_i.
\end{equation}
(\ref{extr3}) follows from (\ref{extr1}) and (\ref{nui}).
\end{proofl}
\begin{remark}\label{rmk2}
Taking $\eta=u-\psi,$ the extremum condition (\ref{extr3}) gives : $\forall i,$
$$\eta u_i^2\lambda_i+(u_i-\psi_i)u_i\nu=0.$$
If the principal curvatures were non-negative, we could conclude : $\forall i,\ u_i^2\leq u_i\psi_i.$ Summing these inequalities and using the Schwarz inequality, we thus would obtain $|\nabla u|\leq |\nabla\psi|.$ Since $|\nabla u|=|Du|\nu$ and $|\nabla\psi|\leq |D\psi|\nu,$ this would thus give the estimate of $\nu$ at the maximum of $\varphi$ :
$$\nu\leq\frac{1}{\sqrt{1-\sup_{\psi\leq\varphi_2}|D\psi|^2}},$$
and thus the estimate of $\nu$ on $B_{R_0}$ since $\eta\geq C_0$ on $B_{R_0},$ where $C_0$ is a positive constant under control. This method, which gives the estimate for the prescribed Gauss curvature equation, was explained to us by O. Schn\"urer. In the prescribed scalar curvature equation, the graph is a priori non-convex. We need to consider a function $\varphi=\eta\nu$ as in (\ref{defeta}), and to use the non-positivity of $\nabla^2\varphi$ at $p$ to estimate further.
\end{remark}
\begin{lemma}\label{lem5}
At $p,$ we have :
\begin{equation}\label{star4}
\left\{\sum_{k}\sigma_{1,k}\lambda_k^2-\sum_k\sigma_{1,k}\lambda_k^2\frac{u_k^2}{\nu^2}\right\}+\sum_i\sigma_{1,i}(\ln\eta)_{ii}\leq|\nabla H|.
\end{equation}
\end{lemma}
\begin{proofl} At $p,$ $\nabla^2\varphi$ is non-positive. We thus have (using summation convention with index ranges $1\leq i,j\leq n):$
$$F^{ij}\varphi_{ij}\leq 0,$$
with $$\varphi_{ij}=\eta_{ij}\nu+\left(\eta_i\nu_j+\eta_j\nu_i\right)+\eta\nu_{ij},$$ and thus
\begin{equation}\label{max}
F^{ij}\nu_{ij}+F^{ij}\frac{\eta_{ij}}{\eta}\nu+F^{ij}\left(\frac{\eta_i}{\eta}\nu_j+\frac{\eta_j}{\eta}\nu_i\right)\leq 0.
\end{equation}
We compute the term $F^{ij}\nu_{ij}$ :
denoting $T=t^k\hat{e_k},$ we have :
\\$\nu_i=-II(T,\hat{e_i})=-t^kh_{ki},$ and $\nu_{ij}=-(t^kh_{ki;j}+t^k_{;j}h_{ki}).$ 
\\Thus
$$F^{ij}\nu_{ij}=-(t^kF^{ij}h_{ki;j}+F^{ij}t^k_{;j}h_{ki}).$$
The Codazzi equations tell us that $h_{ki;j}=h_{ij;k},$ and equation (\ref{eqintr}) differentiated reads : $F^{ij}h_{ij;k}=H_k.$ Thus $F^{ij}h_{ki;j}=H_k.$ Moreover, using the Gauss equation (\ref{gauss2}), we have : $t^k_{;j}=-\nu h_j^k.$  We thus get :
\begin{equation}
F^{ij}\nu_{ij}=-(t^kH_k-\nu F^{ij}h^k_jh_{ki}).
\end{equation}
With the condition (\ref{extr1}),
$$F^{ij}\left(\frac{\eta_i}{\eta}\nu_j+\frac{\eta_j}{\eta}\nu_i\right)=-2F^{ij}\frac{{\eta_i}{\eta_j}}{\eta^2}\nu.$$
We thus obtain the inequality :
\begin{equation}\label{star}
-(t^kH_k-\nu F^{ij}h_j^kh_{ki})+F^{ij}(\ln\eta)_{ij}\nu-F^{ij}\frac{\eta_i\eta_j}{\eta^2}\nu\leq 0,
\end{equation}
Our special choice for the basis $\hat{e_i}$ together with the estimate $t^kH_k\leq\nu|\nabla H|$ and the identity $\displaystyle{\frac{\eta_i}{\eta}=-\frac{u_i\lambda_i}{\nu}}$ thus imply inequality (\ref{star4}).
\end{proofl} 
The next two lemmas are devoted to the estimates of the left hand side terms of (\ref{star4}).

We first obtain an estimate of the second term of (\ref{star4}), by a direct computation using our special choice for the test-function $\eta.$ 

For convenience, let us denote $\xi=u-\psi.$ We thus have $\eta=e^{Ku}(1-e^{-\xi})^K.$
\begin{lemma}\label{lem1}
We have the estimate :
$$\sum_i\sigma_{1,i}(\ln\eta)_{ii}\geq 2KH\nu-\frac{Kc\sigma_1\nu^2\left(2-e^{-\xi}\right)e^{-\xi}}{\left(1-e^{-\xi}\right)^2},$$
where $c$ is a controlled constant.
\end{lemma}
\begin{proofl}
By a direct computation,
$$(\ln\eta)_{ii}=Ku_{ii}+K\left(\frac{\xi_{ii}e^{-\xi}}{1-e^{-\xi}}-\frac{\xi_i^2e^{-\xi}}{(1-e^{-\xi})^2}\right).$$
In view of the Gauss equation  (\ref{gauss2}), $\sum_i\sigma_{1,i}u_{ii}=2H\nu.$ We thus get :
$$\sum_i\sigma_{1,i}(\ln\eta)_{ii}=2KH\nu+K\sum_i\sigma_{1,i}\frac{\xi_{ii}e^{-\xi}}{1-e^{-\xi}}-K\sum_i\sigma_{1,i}\frac{{\xi_i}^2e^{-\xi}}{(1-e^{-\xi})^2}.$$
We readily get the lemma from the following estimates : 
\begin{equation}\label{est}
\sum_i\sigma_{1,i}\xi_i^2\leq c\sigma_1\nu^2, \mbox{ and }|\sum_i\sigma_{1,i}\xi_{ii}|\leq c\sigma_1\nu^2,
\end{equation}
where $c$ is a controlled constant.
\\The inequality $|\nabla\xi|\leq |D\xi|\nu$ implies the first estimate in (\ref{est}). To prove the second estimate in (\ref{est}), let us write 
$$\sum_i\sigma_{1,i}\xi_{ii}=\sum_i\sigma_{1,i}u_{ii}-\sum_i\sigma_{1,i}\psi_{ii}.$$
The first term is $2H\nu$ (Gauss equation). We compute $\psi_{ii}$ following Urbas \cite{U} p. 313. Let $(e_{\alpha})_{1\leq \alpha\leq n}$ denote the standard basis of $\R^{n}.$  We have :
$$\psi_{ii}=\sum_{\alpha}\frac{\partial\psi}{\partial x_{\alpha}}\nabla^2_{ii}(x_{\alpha})+\sum_{\alpha,\beta}\frac{\partial^2\psi}{\partial x_{\alpha}\partial x_{\beta}}\nabla_i(x_{\alpha})\nabla_i(x_{\beta}),$$
and, using the Gauss equations (\ref{gauss1}) and the identity $\nabla_i(x_{\alpha})=<e_{\alpha},\hat{e_{i}}>,$ we obtain :
$$\psi_{ii}=\lambda_i<D\psi,Du>\nu+\sum_{\alpha,\beta}\frac{\partial^2\psi}{\partial x_{\alpha}\partial x_{\beta}}<e_{\alpha},\hat{e_{i}}><e_{\beta},\hat{e_{i}}>.$$
In view of $|<e_{\alpha},\hat{e_{i}}>|\leq\nu,$ this gives the second estimate of (\ref{est}), and the lemma \ref{lem1}.
\end{proofl}

We next estimate the first term of the left hand side term of (\ref{star4}), in the spirit of the method used to estimate the gradient for the prescribed mean curvature equation in lorentzian manifolds (see for instance \cite{B},\cite{CY}). For that purpose, we use the following inequality, which is a particular case of an inequality of Ivochkina, Lin and Trudinger; it is proved in \cite{LT} (ineq. (26)): if $\lambda_{i_0}\leq 0,$ then
\begin{equation}\label{ineqLT}
\varepsilon \sigma_{1,i_0}\lambda_{i_0}^2\leq\sum_{k\neq i_0}\sigma_{1,k}\lambda_k^2,
\end{equation}
where $\varepsilon=\varepsilon(n)$ is a small constant. We deduce the following inequality : 
\begin{lemma}\label{lem4} If $\lambda_{i_0}\leq 0,$
\begin{equation}\label{csqineqLT}
\sum_{k}\sigma_{1,k}\lambda_k^2-\sum_k\sigma_{1,k}\lambda_k^2\frac{u_k^2}{\nu^2}\geq \varepsilon\sigma_{1,i_0}\lambda_{i_0}^2.
\end{equation}
\end{lemma}
\begin{proofl}
Let us write :
\begin{equation}\label{eqn1lem4}
\sum_k\sigma_{1,k}\lambda_k^2\frac{u_k^2}{\nu^2}=\sum_{k\neq i_0}\sigma_{1,k}\lambda_k^2\frac{u_k^2}{\nu^2}+ \sigma_{1,i_0}\lambda_{i_0}^2\frac{u_{i_0}^2}{\nu^2}.
\end{equation}
The first term in (\ref{eqn1lem4}) is smaller than $\displaystyle{\left(\sum_{k\neq i_0}\sigma_{1,k}\lambda_k^2\right)\left(\sum_{k\neq i_0}\frac{u_k^2}{\nu^2}\right).}$
\\Using (\ref{ineqLT}), the second term in (\ref{eqn1lem4}) is smaller than 
$$(1-\varepsilon)\sigma_{1,i_0}\lambda_{i_0}^2\frac{u_{i_0}^2}{\nu^2}+\left(\sum_{k\neq{i_0}}\sigma_{1,k}\lambda_k^2\right)\frac{u_{i_0}^2}{\nu^2}.$$
We thus have 
$$\sum_k\sigma_{1,k}\lambda_k^2\frac{u_k^2}{\nu^2}\leq \left(\sum_{k\neq i_0}\sigma_{1,k}\lambda_k^2\right)\left(\sum_{k\neq i_0}\frac{u_k^2}{\nu^2}+\frac{u_{i_0}^2}{\nu^2}\right)+(1-\varepsilon)\sigma_{1,i_0}\lambda_{i_0}^2\frac{u_{i_0}^2}{\nu^2},$$
and thus
$$\sum_k\sigma_{1,k}\lambda_k^2\frac{u_k^2}{\nu^2}\leq \sum_{k\neq i_0}\sigma_{1,k}\lambda_k^2+(1-\varepsilon)\sigma_{1,i_0}\lambda_{i_0}^2\frac{u_{i_0}^2}{\nu^2}.$$
Thus
$$\sum_{k}\sigma_{1,k}\lambda_k^2-\sum_k\sigma_{1,k}\lambda_k^2\frac{u_k^2}{\nu^2}\geq \sigma_{1,i_0}\lambda_{i_0}^2\left(1-(1-\varepsilon)\frac{u_{i_0}^2}{\nu^2}\right),$$
which gives the result since $u_{i_0}^2\leq \nu^2.$
\end{proofl}
Multiplying (\ref{star4}) by $(1-e^{-\xi})^2$ and using the lemmas \ref{lem1} and \ref{lem4}, we thus finally obtain the inequality : if $\lambda_{i_0}\leq 0,$
\begin{equation}\label{star2}
\begin{split}\varepsilon\sigma_{1,i_0}\lambda_{i_0}^2(1-e^{-\xi})^2+2KH\nu(1-e^{-\xi})^2-Kc\sigma_1\nu^2(2-e^{-\xi})e^{-\xi}\\\leq|\nabla H|(1-e^{-\xi})^2.
\end{split}
\end{equation}
With this inequality at hand, we can prove lemma \ref{lem} and thus obtain the gradient estimate.
\\
\\{\textit Proof of lemma }\ref{lem} : 
Let $\alpha>0$ such that $1-\alpha>\sup_{\psi\leq\varphi_2}|D\psi|.$ We first prove by contradiction that for $K$ sufficiently large, we have at $p,$ maximum interior of $\varphi$ : $\forall i,$ 
\begin{equation}
(1-\alpha)u_i^2\leq\psi_i u_i.\label{lastineq}
\end{equation}
We thus suppose that there exists $i_0$ such that $(1-\alpha)u_{i_0}^2>\psi_{i_0}u_{i_0},$ ie such that $\frac{\xi_{i_0}}{u_{i_0}}>\alpha.$ We recall that $\eta=e^{Ku}\left(1-e^{-\xi}\right)^K.$ Lemma \ref{lem0} for $i_0$ gives
\begin{eqnarray}
\lambda_{i_0} & = & -K\nu\left(1+\frac{\xi_{i_0}}{u_{i_0}}\frac{e^{-\xi}}{1-e^{-\xi}}\right)\nonumber\\
& \leq & -K\nu\alpha\frac{e^{-\xi}}{1-e^{-\xi}}.\nonumber
\end{eqnarray}
Thus $\lambda_{i_0}\leq 0,$ and $\lambda_{i_0}^2\geq K^2\nu^2\alpha^2\frac{e^{-2\xi}}{(1-e^{-\xi})^2}.$
\\
\\Inequality (\ref{star2}) thus gives
\begin{equation}\label{star3}
\begin{split}
\varepsilon\sigma_{1,i_0}K^2\nu^2\alpha^2e^{-2\xi}+2K\inf_{\Gamma}H\nu(1-e^{-\xi})^2-Kc\sigma_1\nu^2(2-e^{-\xi})e^{-\xi}\\\leq \sup_{\Gamma}|\nabla H|(1-e^{-\xi})^2.\end{split}
\end{equation}
Moreover, we have the estimates :
$$\inf_{\Gamma}H\geq \inf_{\psi\leq\varphi_2}H\mbox{ and }\sup_{\Gamma}|\nabla H|\leq \nu\sup_{\psi\leq\varphi_2}|DH|_{eucl}.$$
Since $\lambda_{i_0}\leq 0,$ we have $\sigma_{1,i_0}\geq \sigma_1,$ and (\ref{star3}) is impossible for $K$ sufficiently large. This achieves the proof of (\ref{lastineq}). We then achieve the proof of the lemma : as in the convex case (see remark \ref{rmk2} above), summing these inequalities and using the Schwarz inequality, we get  $(1-\alpha)|\nabla u|\leq |\nabla\psi|.$ Since $|\nabla u|=\sqrt{\nu^2-1}$ and $|\nabla\psi|\leq |D\psi|\nu,$ this gives the estimate (\ref{estnu}) of $\nu$ at the maximum of $\varphi.$ 
\begin{remark}
 There is an alternative presentation of this estimate if we restrict our study to the case of hypersurfaces asymptotic to the lightcone at infinity, using the other time function $\tau=\sqrt{{x_{n+1}}^2-|x|^2}.$ The proof is more elegant, and some technical difficulties disappear.  See Bartnik \cite{B2} for the mean curvature case. Nevertheless, we need here a more general estimate allowing us to study asymptotic values which are bounded perturbations of the lightcone at infinity (see section 6). To this purpose, the presentation above is convenient. We hope to come back soon to these questions. 
\end{remark}

\section{The uniform $C^2$ estimate}\label{c2}
Our uniform $C^2$ estimate relies on the work of Urbas \cite{U} where the author obtains a $C^2$ interior estimate for solutions of the Dirichlet problem for the prescribed scalar curvature equation in $\R^{n,1},$ in the case of a spacelike affine boundary data; its proof uses in a crucial way the uniform convexity of an extension of the boundary data.

Following Urbas proof, we easily prove the $C^2$ estimate :
\begin{proposition}\label{propu}
There exist $C_{R_0}$ and $R_1=R_1(R_0)$ such that $\forall R\geq R_1,$
$$\|u_R\|_{2,B_{R_0}}\leq C_{R_0}.$$
\end{proposition}
\begin{proofp} Let $\psi=-1+\sqrt{a^2+|x|^2},$ with $a=a(R_0)$ sufficiently large such that $\psi\geq\varphi_2+1$ on $B_{R_0}.$ Since $\psi$ is asymptotic to the cone $x_{n+1}=-1+|x|,$ we have $\psi<\varphi_1$ out a large fixed ball. For every $R$ the set $\{\psi\geq u_R\}$ is in this fixed ball, and thus, from the gradient estimate, $|Du_R|$ is under control on the set $\{\psi\geq u_R\}$ if $R$ is large, ie :
$$\exists\ \theta_{R_0}\ \in\ [0,1[\ and\ R_1=R_1(R_0)\ s.t.\ \sup_{\{\psi\geq u_R\}}|Du_R|\leq 1-\theta_{R_0}\ if\ R\geq R_1.$$ 
Let $R\geq R_1$ and let us denote $u_R$ by $u.$ Let $\Gamma$ be the piece of the graph $M$ of $u$ on which $\psi\geq u,$ and set $\eta=\psi-u,$ considered as a function on $\Gamma.$ Following Urbas \cite{U} without any modification, the function $\displaystyle{W(x,\xi)=\eta^{\beta}h_{ij}\xi^i\xi^j}$ defined for all $x\ \in\ \Gamma$ and all unit $\xi\ \in\ T_xM$ is bounded by a controlled constant at $(x_0,\xi_0)$ where it attains its maximum, if $\beta$ is a large constant (here, it is crucial that $\psi$ be uniformly convex on $\{\psi\geq\varphi_1\}$ and thus on $\{\psi\geq u\}$; see \cite{U}). Since $\eta=\psi-u$ is larger than 1 above $B_{R_0},$ this gives a $C^2$ estimate of $u$ on $B_{R_0}.$
\end{proofp}
\section{General data at infinity}\label{lastsec}
Let us turn to the proof of theorem \ref{the1}. The first paragraph is devoted to the construction of barriers, following Treibergs \cite{T}. The proof then reduces to obtain uniform gradient and $C^2$ estimates. For that purpose, following the proofs in sections \ref{gradient} and \ref{c2}, we only need to construct convenient auxiliary functions. It is the aim of the second and the third paragraph.
\subsection{Construction of barriers}
\ 
\\We follow closely the Treibergs construction \cite{T}. We assume that the function $f$ is not constant, and we extend it to $\R^n\backslash\{0\}$ by $f(x)=f\left(\frac{x}{|x|}\right).$ Since $f$ is $C^2,$ there exists $M$ such that : $\forall x,y\in S^{n-1},$
$$|f(x)-f(y)-Df(y).(x-y)|\leq M|x-y|^2=-2My.(x-y).$$
Setting, for $i=1,2,$
$$p_i(y)=Df(y)+(-1)^{i+1}2My,$$
we get : $\forall x,y\in S^{n-1},$
$$p_1(y).(x-y)\leq f(x)-f(y)\leq p_2(y).(x-y).$$
We set, for $i=1,2,$ and $\forall x\in\R^n,$ $\forall y\in S^{n-1},$
$$z_i(x,y)=f(y)-p_i(y).y+\left(\alpha_i^2+|x+p_i(y)|^2\right)^{\frac{1}{2}},$$
where $\alpha_i=\sqrt{\frac{n(n-1)}{2h_i}}$ (see remark 1.3 above). It is easy to check that, $\forall x,y\in S^{n-1},$
$$\lim_{r\rightarrow+\infty}z_1(rx,y)-r\leq f(x)\leq \lim_{r\rightarrow+\infty}z_2(rx,y)-r,$$
where the inequalities become equalities when $x=y.$ We then set, $\forall x\in\R^n :$ 
$$q_1(x)=\sup_{y\in S^{n-1}}z_1(x,y),\ q_2(x)=\inf_{y\in S^{n-1}}z_2(x,y).$$
We have : $\forall x\in \R^n,$
$$q_1(x)<q_2(x),$$
and, for $i=1,2,$
\begin{equation}\label{limq}
\lim_{|x|\rightarrow +\infty}q_i(x)-|x|-f\left(\frac{x}{|x|}\right)=0.
\end{equation}
The functions $q_1$ and $q_2$ are continuous, but do not belong to $C^1(\R^n)$ in general. As easily seen, $q_1$ is a strictly convex function, and for $m\in\N$ (sufficiently large), there exists $\Omega_m,$ open set $C^{\infty}$, strictly convex, such that
$$\{x|\ q_2(x)<m\}\subset\Omega_m\subset\{x|\ q_1(x)<m\}.$$
We then solve the Dirichlet problems $\scal[u_m]=H(.,u_m)$ in $\Omega_m,$ and $u_m=m$ on $\partial\Omega_m,$ relying on theorem 1.1 of Urbas \cite{U} and a classical use of the Leray-Schauder fixed point theorem. Specifically, let $\varphi_m$ be a smooth and strictly convex spacelike function on $\Omega_m,$ equal to $m$ on $\partial\Omega_m.$ Let us consider the Banach space
$$E=\{v\in\ C^{2,\alpha}(\overline{\Omega_m})|\ v=0\mbox{ on }\partial\Omega_m\}.$$
Thanks to theorem 1.1 of \cite{U}, we may define the operator
$$T:E\times[0,1]\rightarrow E,\ (v,\sigma)\mapsto w,$$
where $w$ is the function $C^{4,\alpha}$ such that $w+\varphi_m$ is admissible and satisfies :
$$\scal[w+\varphi_m]=\sigma H(.,v+\varphi_m)+(1-\sigma)\scal[\varphi_m]\mbox{ in }\Omega_m,\mbox{ and }w=0\mbox{ on }\partial\Omega_m.$$
$T$ is a compact map (thanks to the gradient and the $C^{2,\alpha}$ a priori estimates), $T(.,0)$ is the null operator, and the fixed points $w$ solutions of $T(w,\sigma)=w$ admit a $C^{2,\alpha}(\overline{\Omega_m})$ bound independent of $\sigma\in\ [0,1]$ (see the estimates obtained in the proof of theorem \ref{theo3} section \ref{sec00} above). Thus, the Leray-Schauder theorem yields a fixed point $w$ for $T(.,1).$ The function $u_m=w+\varphi_m$ is then a solution of the Dirichlet problem. 
\\On $\partial\Omega_m,$ we have $q_1\leq m$ and $q_2\geq m;$ we thus get : $\forall y\ \in S^{n-1},$
$$z_1(.,y)\leq u_m\leq z_2(.,y)\mbox{ on }\partial\Omega_m.$$
Since the graphs of $z_1(.,y)$ and $z_2(.,y)$ are hyperbolo\"{\i}ds whose $\sigma_2$ are respectively equal to $h_1$ and $h_2,$ the pinching condition on $H$ readily implies : $\forall y\ \in S^{n-1},$
$$\scal[z_1(.,y)]\geq\scal[u_m]\geq\scal[z_2(.,y)]\mbox{ on }\Omega_m.$$
The comparison principle applied to the operator $\scal$ implies : $\forall y\ \in S^{n-1},$
$$z_1(.,y)\leq u_m\leq z_2(.,y)\mbox{ on }\Omega_m,$$
and we thus obtain : $\forall m,$
$$q_1\leq u_m\leq q_2.$$
To prove theorem \ref{the1}, it only remains to construct uniform estimates for the gradient and for the second derivatives on every compact subset.
\subsection{The uniform gradient estimate}
\ 
\\We apply the procedure of section \ref{gradient}. To that purpose, we only have to construct a function $\psi$ satisfying properties as in (\ref{ptepsi}). It is the aim of the following lemma :
\begin{lemma}
Let $B_{R_0}$ be a fixed ball of radius $R_0.$ There exist $\delta>0,$ $R_1,R_2$ with $R_2>R_1>R_0,$ and $\psi$ spacelike on $B_{R_2}$ such that :
$$\psi\leq q_1-\delta\mbox{ on }B_{R_0},$$
and
$$\psi>q_2\mbox{ on }B_{R_2}\backslash B_{R_1}.$$ 
\end{lemma}
\begin{proofl}
We first construct $\delta>0,$ independent of $y\in S^{n-1},$ such that the difference of heights between the hyperbolo\"{\i}d $x_{n+1}=z_1(x,y)$ and its asymptotic cone is larger than $\delta$ on $B_{R_0}:$ we remark to this purpose that the hyperbolo\"{\i}d $x_{n+1}=z_1(x,y)$ is translated from the hyperbolo\"{\i}d $x_{n+1}=\sqrt{\alpha_1^2+|x|^2}$ by the vector $(-p_1(y),f(y)-p_1(y).y).$ Let $R$ be sufficiently large such that $\forall y\in S^{n-1},$
$$B_{R_0}\subset \overline{B}(-p_1(y),R).$$
Setting 
$$\delta_0=\inf_{B_R}\left\{\left(\alpha_1^2+|x|^2\right)^{\frac{1}{2}}-|x|\right\},$$
every $\delta<\delta_0$ is suitable.
\\
\\We then set
$$z_1^{\varepsilon}(x,y)=f(y)+\varepsilon-p_1(y).y+\left(\varepsilon^2+|x+p_1(y)|^2\right)^{\frac{1}{2}}.$$
Taking $\varepsilon$ small such that
$$\varepsilon+\left(\varepsilon^2+|x|^2\right)^{\frac{1}{2}}<\left(\alpha_1^2+|x|^2\right)^{\frac{1}{2}}-\delta\mbox{ on }\overline{B}(0,R),$$
we see that :
$$z_1^{\varepsilon}<z_1-\delta\mbox{ on }B_{R_0}.$$
Setting
$$q_1^{\varepsilon}(x)=\sup_{y\in S^{n-1}}z_1^{\varepsilon}(x,y),$$
we thus have :
\begin{equation}\label{ineq1}
q_1^{\varepsilon}<q_1-\delta\mbox{ on }B_{R_0},
\end{equation}
and, since we get in view of (\ref{limq})
$$\lim_{|x|\rightarrow +\infty}q_1^{\varepsilon}(x)-q_2(x)=\varepsilon,$$
we also have :
\begin{equation}\label{ineq2}
q_1^{\varepsilon}\geq q_2+\frac{\varepsilon}{2}\mbox{ on }\R^n\backslash B_{R_1},
\end{equation}
for $R_1$ sufficiently large.
\\
\\We finally need to regularize $q_1^{\varepsilon}.$ We set $R_2>R_1.$ We first note that there exists $\theta\in\ ]0,1[$ such that $q_1^{\varepsilon}$ is lipschitz continuous with coefficient $1-\theta$ on $B_{R_2+1}.$ Let $\varphi\in C^{\infty}_c(\R^n)$ be such that $\varphi\geq 0,$ and $\int_{\R^n}\varphi=1,$ and set $\varphi_{\beta}(x)=\frac{1}{\beta^n}\varphi(\frac{x}{\beta}).$ The convolution $q_1^{\varepsilon}\star \varphi_{\beta}$ is smooth, lipschitz continuous with coefficient $1-\theta$ on $B_{R_2}$ (and thus spacelike), and converges to $q_1^{\varepsilon}$ uniformly on $B_{R_2}$ when $\beta$ tends to 0. In view of (\ref{ineq1}) and (\ref{ineq2}), for $\beta$ sufficiently small, $\psi=q_1^{\varepsilon}\star \varphi_{\beta}$ satisfies the conditions of the lemma.
\end{proofl}
\subsection{The uniform $C^2$ estimate}
\ \\
We use the work of Urbas \cite{U}, as in section \ref{c2}. For that purpose, recording the proof of proposition \ref{propu}, we only need the following lemma :
\begin{lemma}
There exist $R_1,R_2$ with $R_2>R_1>R_0$ and $\psi,$ strictly convex, spacelike on $B_{R_2},$ such that 
$$\psi\geq q_2+1\mbox{ on }B_{R_0},$$
and
$$\psi<q_1\mbox{ on }B_{R_2}\backslash B_{R_1}.$$
\end{lemma}
\begin{proofl}
We set $M=\sup_{B_{R_0}}q_2+2,$ and we fix $\varepsilon,$ a positive constant.
\\We define
$$z_1^{-\varepsilon}(x,y)=f(y)-\varepsilon-p_1(y).y+\left(\alpha^2+|x+p_1(y)|^2\right)^{\frac{1}{2}},$$
with $\alpha$ sufficiently large such that $\forall y\in S^{n-1},$
$$z_1^{-\varepsilon}(x,y)\geq M\mbox{ on }B_{R_0}.$$
Setting 
$$q_1^{-\varepsilon}(x)=\sup_{y\in S^{n-1}}z_1^{-\varepsilon}(x,y),$$
we have $q_1^{-\varepsilon}\geq M$ on $B_{R_0},$ and (see (\ref{limq}))
$$\lim_{|x|\rightarrow +\infty}q_1^{-\varepsilon}(x)-q_1(x)=-\varepsilon.$$
Thus $\displaystyle{q_1^{-\varepsilon}(x)\leq q_1(x)-\frac{\varepsilon}{2}}$ for all $x\in\R^n\backslash B_{R_1},$ if $R_1$ is chosen sufficiently large.
\\We regularize $q_1^{-\varepsilon}$ as in the previous section. We choose $R_2>R_1.$ Taking $\varphi_{\beta}$ as above, the convolution $q_1^{-\varepsilon}\star\varphi_{\beta}$ is smooth, spacelike, strictly convex, and is uniformly convergent to $q_1^{-\varepsilon}$ on $B_{R_2}.$ Thus, if $\beta$ is sufficiently small, $\psi=q_1^{-\varepsilon}\star\varphi_{\beta}$ satisfies the conditions of the lemma.
\end{proofl}
\\
\\
\textit{Acknowledgements : }I would like to thank O. Schn\"urer for numerous discussions which are at the origin of this work. The uniform interior gradient estimate required to prove existence of entire hypersurfaces of prescribed Gauss curvature in Minkowski space which are asymptotic to a lightcone is due to him. I adapted it here to the scalar curvature case (see the introduction and remark \ref{rmk2} above). Ph. Delano\"e mentioned to me in 2001 the problem of finding entire solutions of prescribed curvature as a natural continuation of my doctoral dissertation. I thank him for this, and for numerous useful comments which helped me to considerably improve the results and the writing of this article. 


\begin{thebibliography}{99}
\bibitem[1] {B} R.Bartnik. Existence of maximal surfaces in asymptotically flat spacetimes. Commun. Math. Phys. \textbf{94} (1984) 155-175.
\bibitem[2] {B2} R.Bartnik. Maximal surfaces and general relativity, Geometry and partial differential equations, 2nd Miniconf., Canberra/Aust. 1986, Proc. Cent. Math. Anal. Aust. Natl. Univ. \textbf{12} (1987) 24-49. 
\bibitem[3] {Bay} P.Bayard. Dirichlet problem for space-like hypersurfaces with prescribed scalar curvature in $\R^{n,1}$. Calc. Var. \textbf{18} (2003) 1-30.
\bibitem[4] {CH} R.Courant, D.Hilbert. Methods of mathematical physics, Vol. 2. Wiley Interscience Publishers (1962).
\bibitem[5] {CY} S.Y.Cheng, S.T.Yau. Maximal space-like hypersurfaces in the Lorentz-Minkowski spaces. Ann. Math.(2) \textbf{104} (1976) 407-419.
\bibitem[6] {De} Ph.Delano\"e. Dirichlet problem for the equation of prescribed Lorentz-Gauss curvature. Ukr. Math. J. \textbf{42:12} (1990) 1538-1545.
\bibitem[7] {De2} Ph.Delano\"e. Personal communication. March 2005.
\bibitem[8] {G} B.Guan. The Dirichlet problem for Monge-Amp\`ere equations in non-convex domains and spacelike hypersurfaces of constant Gauss curvature. Trans. Amer. Math. Soc. \textbf{350} (1998) 4955-4971.
\bibitem[9] {GT} D.Gilbarg, N.S.Trudinger. Elliptic partial differential equations of second order, second edition, revised third printing, Springer-Verlag (1998).
\bibitem[10] {L} A.M.Li. Spacelike hypersurfaces with constant Gauss-Kronecker curvature in the Minkowski space. Arch. Math. \textbf{64:6} (1995) 534-551.
\bibitem[11] {LT} M.Lin, N.S.Trudinger. On some inequalities for elementary symmetric functions. Bull. Austral. Math. Soc. \textbf{50} (1994) 317-326.
\bibitem[12] {BP} A.J.B.Potter, An elementary version of the Leray-Schauder theorem. J. London Math. Soc. \textbf{5:2} (1972) 414-416.  
\bibitem[13] {Sch} O.Schn\"urer. The Dirichlet problem for Weingarten hypersurfaces in Lorentzian manifolds. Math. Z. \textbf{242} (2002) 159-181.
\bibitem[14] {T} A.Treibergs. Entire Spacelike hypersurfaces of constant mean curvature in Minkowski space. Invent. Math. \textbf{66} (1982) 39-56.
\bibitem[15] {U} J.Urbas. The Dirichlet problem for the equation of prescribed scalar curvature in Minkowski space. Calc. Var. \textbf{18} (2003) 307-316.
\end{thebibliography}
\end{document}